\documentclass[12pt]{amsart}

\usepackage{amssymb}
\usepackage[curve]{xypic}
\usepackage{amsmath}
\usepackage{amsthm}
\usepackage{amsfonts}
\usepackage{amscd}
\usepackage{pifont}

\newtheorem{Lemma}{Lemma}[section]
\newtheorem{Theorem}[Lemma]{Theorem}
\newtheorem{Corollary}[Lemma]{Corollary}
\newtheorem{Proposition}[Lemma]{Proposition}
\newtheorem*{LemmaAppendix1}{Lemma \ref{LRB: dimensions of Ext spaces}}

\theoremstyle{definition}

\newtheorem{Remark}[Lemma]{Remark}
\newtheorem{Example}[Lemma]{Example}

\numberwithin{equation}{section}

\newcommand{\supp}{\operatorname{supp}}

\newcommand{\Hom}{\operatorname{Hom}}
\newcommand{\Ext}{\operatorname{Ext}}

\renewcommand{\AA}{\mathcal A}
\newcommand{\FF}{\mathcal F}
\newcommand{\LL}{\mathcal L}

\newcommand{\tot}[2]{\buildrel #1 \over #2}

\allowdisplaybreaks

\begin{document}

\title{The Quiver of the Semigroup Algebra of a Left Regular Band}
\author{Franco V Saliola}
\email{saliola@gmail.com}

\begin{abstract}
Recently it has been noticed that many interesting 
combinatorial objects belong to a class of semigroups called \emph{left
regular bands}, and that random walks on these semigroups encode several
well-known random walks. For example, the set of faces of a hyperplane
arrangement is endowed with a left regular band structure. This paper
studies the module structure of the semigroup algebra of an arbitrary left
regular band, extending results for the semigroup algebra of the faces of a
hyperplane arrangement. In particular, a description of the quiver of the
semigroup algebra is given and the Cartan invariants are computed. These
are used to compute the quiver of the face semigroup algebra of a
hyperplane arrangement and to show that the semigroup algebra of the free
left regular band is isomorphic to the path algebra of its quiver.
\end{abstract}

\maketitle

\markboth{\textsc{FRANCO V SALIOLA}}{\textsc{THE SEMIGROUP ALGEBRA OF A
LEFT REGULAR BAND}}

\small
\tableofcontents
\normalsize

\section{Introduction}
A \emph{left regular band} is a semigroup $S$ satisfying $x^2 = x$ and $xyx
= xy$ for all $x,y \in S$. Recent interest in left regular bands and their
semigroup algebras arose due to the work of K. S. Brown \cite{Brown2000},
in which the representation theory of the semigroup algebra is used to
study random walks on the semigroup. There are several interesting examples
of such random walks, including the random walk on the chambers of a
hyperplane arrangement. Several detailed examples are included in
\cite{Brown2000}.

The starting point of this paper is the fact that the irreducible
representations of the semigroup algebra of a left regular band are all
1-dimensional. This implies that there is a canonical quiver (a directed
graph) associated to the left regular band, and that the semigroup algbera
is a quotient of the path algebra of the quiver. This paper determines a
combinatorial description of this quiver and the Cartan invariants of the
semigroup algebras and illustrates the theory through detailed examples.

The paper is structured as follows. Section \ref{section: LRBs} recalls the
definition and collects some properties of left regular bands, and
introduces the examples that will be used throughout the paper. Section
\ref{section: representations of the semigroup algebra} describes the
irreducible representations of the semigroup algebra of a left regular
band. In Section \ref{LRB: section: primitive idempotents} a complete
system of primitive orthogonal idempotents for the semigroup algebra is
explicitly constructed. Section \ref{LRB: section: projective
indecomposable modules} describes the projective indecomposable modules of
the semigroup algebra. Sections \ref{LRB: quiver} through \ref{section:
example: the face semigroup of a hyperplane arrangement} deal with
computing the quiver of the semigroup algebra. Sections \ref{section:
idempotents in the subalgebras} through \ref{section: the cartan invariants
for a free left regular band} compute the Cartan invariants of the
semigroup algebras.  Finally, Section \ref{section: future directions}
discusses future directions for this project.

\section{Left Regular Bands}
\label{section: LRBs}

See \cite[Appendix B]{Brown2000} for foundations of left regular bands and for
proofs of the statements presented in this section.

A \emph{left regular band} is a semigroup $S$ satisfying the following
two properties.
\begin{quote}
\begin{itemize}
\item[(LRB1)] $x^2 = x$ for all $x \in S$.
\item[(LRB2)] $xyx = xy$ for all $x,y \in S$.
\end{itemize}
\end{quote}

Define a relation on the elements of $S$ by $y \leq x$ iff $yx = x$.
This relation is a partial order (reflexive, transitive and antisymmetric),
so $S$ is a poset.

Define another relation on the elements of $S$ by $y \preceq x$ iff
$xy = x$. This relation is reflexive and transitive, but not necessarily
antisymmetric. Therefore we get a poset $L$ by identifying $x$ and
$y$ if $x \preceq y$ and $y \preceq x$. Let $\supp: S \to L$ denote the
quotient map. $L$ is called the \emph{support semilattice} of $S$ and
$\supp: S \to L$ is called the \emph{support map}.

\begin{Proposition}
\label{LRB: associated lattice}
If $S$ is a left regular band, then there is a semilattice $L$ and
a surjection $\supp: S \to L$ satisfying the following properties for
all $x,y \in S$.
\begin{enumerate}
 \setlength{\itemsep}{1pt} \setlength{\parskip}{1pt} \setlength{\parsep}{1pt}
\item If $y \leq x$, then $\supp(y) \leq \supp(x)$.
\item $\supp(xy) = \supp(x) \vee \supp(y)$.
\item 
\label{LRB: xy=x iff XleqY}
$xy = x$ iff $\supp(y) \leq \supp(x)$.
\item If $S'$ is a subsemigroup of $S$, then the image of $S'$ in $L$
 is the support semilattice of $S'$.
\end{enumerate}
\end{Proposition}
Statement (1) says that $\supp$ is an order-preserving poset map. (2) says
that $\supp$ is a semigroup map where we view $L$ as a semigroup with
product $\vee$. (3) follows from the construction of $L$, and (4) follows
from the fact that (3) characterizes $L$ upto isomorphism.
If $S$ has an identity element then $L$ has a minimal
element $\hat 0$. If, in addition, $L$ is finite, then $L$ has a maximal
element $\hat 1$, and is therefore a lattice
\cite[Proposition 3.3.1]{Stanley1997}. In this case $L$ is the
\emph{support lattice} of $S$.

\begin{Example}[The Free Left Regular Band]
\label{LRB: free LRB}
The \emph{free left regular band} $F(A)$ \emph{with identity} on
a finite set $A$ is the set of all (ordered) sequences of distinct
elements from $A$ with multiplication defined by
\begin{align*}
(a_1, \ldots, a_l) \cdot (b_1, \ldots, b_m) =
(a_1, \ldots, a_l, b_1, \ldots, b_m)^\text{\ding{36}}
\end{align*}
where \ding{36} means ``delete any element that has occured earlier''.
Equivalently, $F(A)$ is the set of all words on the alphabet $A$ that
do not contain any repeated letters. 

\begin{figure}
\begin{gather*}
\xymatrix@C=0pt@R=1em{
abc \ar@{-}[d] && acb \ar@{-}[d] && bac \ar@{-}[d] && bca \ar@{-}[d] & & cab \ar@{-}[d] && cba \ar@{-}[d] \\
 ab \ar@{-}[dr] &&  ac \ar@{-}[dl] &&  ba \ar@{-}[dr] &&  bc \ar@{-}[dl] & &  ca \ar@{-}[dr] &&  cb \ar@{-}[dl] \\
    &a \ar@{-}[drrrr] &    &&    &b\ar@{-}[d]&     & &    &c\ar@{-}[dllll]& \\
    &&     &&    &1& 
} 
\hspace{3em}
\xymatrix@C=1em@R=1em{
  & abc \ar@{-}[dl] \ar@{-}[rd] \ar@{-}[d] \\
ab \ar@{-}[d]\ar@{-}[dr] & ac \ar@{-}[dl]\ar@{-}[dr] & bc \ar@{-}[d]\ar@{-}[dl] \\
a \ar@{-}[dr]  & b \ar@{-}[d]  & c \ar@{-}[dl] \\
 & \emptyset
}
\end{gather*}
\caption[The free left regular band on three generators]
{The poset of the free left regular band $F(\{a,b,c\})$ on
three generators and its support lattice.}
\label{figure: free lrb}
\end{figure}
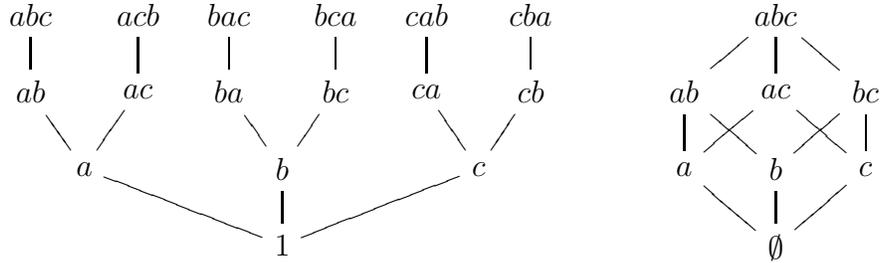

The empty sequence is an element of $F(A)$, therefore
$F(A)$ contains an identity element.
The support lattice of $F(A)$ is the lattice $L$ of subsets of $A$
and the support map $\supp: F(A) \to A$ sends a sequence $(a_1, \ldots,
a_l)$ to the set of elements in the sequence $\{a_1, \ldots, a_l\}$.
Figure \ref{figure: free lrb} shows the Hasse diagrams of the poset
$(F(A),\leq)$ and the support lattice of $F(A)$, where $A = \{a,b,c\}$.
\end{Example}

\begin{Example}[Hyperplane Arrangements]
\label{LRB: hyperplane arrangements}
A \emph{(central) hyperplane arrangement} $\AA$ is a finite collection of
hyperplanes containing the origin in some real vector space $V = \mathbb
R^d$, for some $d \in \mathbb N$.
For each hyperplane $H \in \AA$, let $H^+$ and $H^-$ denote the two open
half spaces of $V$ determined by $H$. The choice of labels $H^+$ and $H^-$
on the two open half spaces is arbitrary, but fixed throughout. For
convenience, let $H^0$ denote $H$. A \emph{face} of the arrangement $\AA$
is a non-empty intersection of the form $\bigcap_{H \in \AA}
H^{\epsilon_H}$, where $\epsilon_H \in \{0,+,-\}$. Let $\FF$ denote the set
of all faces of $\AA$. Define a relation on $\FF$ by $x \leq y$ iff $x
\subset \overline{y}$, where $\overline{y}$ denotes the closure of the set
$y$. The relation is a partial order.

If $x = \bigcap_{H \in \AA} H^{\epsilon_H}$ is a face, then let
$\sigma_H(x) = \epsilon_H$ and let $\sigma(x) = (\sigma_H(x))_{H \in \AA}$.
The sequence $\sigma(x)$ is called the \emph{sign sequence} of $x$.
Define the product of two faces $x, y \in \FF$ to be the face $xy$
with sign sequence
\begin{gather*}
\sigma_H(xy) = 
\begin{cases}
\sigma_H(x), & \text{if } \sigma_H(x) \neq 0, \\
\sigma_H(y), & \text{if } \sigma_H(x) = 0.
\end{cases}
\end{gather*}
This product has a geometric interpretation: the product $xy$ of two faces
$x,y$ is the face entered by moving a small positive distance along a
straight line from any point in $x$ to a point in $y$.  It is
straightfoward to verify that this product gives $\FF$ the structure of an
associative left regular band.
Since all the hyperplanes in the arrangement contain the origin,
$\FF$ contains an identity element: $\cap_{H \in \AA} H$.
The left regular band $\FF$ is called the \emph{face semigroup} of $\AA$,
and the semigroup algebra $k\FF$ of $\FF$ is called the face semigroup
algebra of $\AA$.

Let $\LL$ denote the set of subspaces of $V$ that can be obtained as the
intersection of some hyperplanes in $\AA$. Then $\LL$ is a finite lattice,
called the \emph{intersection lattice} of $\AA$, where the subspaces are
ordered by inclusion and the meet operation is intersection. (Note that
some authors order $\LL$ by reverse inclusion rather than inclusion.)
$\LL$ is the support lattice of $\FF$ and
the support map $\supp: \FF \to \LL$ maps
a face $x \in \FF$ to the intersection of all the hyperplanes of the
arrangement that contain the face:
$\supp(x) = \bigcap_{\{H \in \AA: x \subset H\}} H$.
\end{Example}

\section{Representations of the Semigroup Algebra}
\label{section: representations of the semigroup algebra}

Let $k$ denote a field and $S$ a left regular band. The semigroup algebra
of $S$ is denoted by $kS$ and consists of all formal linear combinations
$\sum_{s \in S} \lambda_s s$, with $\lambda_s \in k$ and multiplication
induced by $\lambda_s s \cdot \lambda_t t = \lambda_s \lambda_t st$, where
$st$ is the product of $s$ and $t$ in the semigroup $S$.
The following summarizes Section 7.2 of \cite{Brown2000}.

Since $S$ and $L$ are semigroups and $\supp: S \to L$ is a semigroup
morphism, the support map 
extends linearly to a surjection of semigroup algebras $\supp: kS \to kL$. 
The kernel of this map is nilpotent and the semigroup algebra
$kL$ is isomorphic to a product of copies of the field $k$,
one copy for each element of $L$. 
Standard ring theory implies that $\ker(\supp)$ is the Jacobson radical of $kS$
and that 
the irreducible
representations of $kS$ are given by the components of the composition
$kS \tot{\supp}\longrightarrow kL \tot{\cong}\longrightarrow 
\prod_{X\in L} k$. This last map sends $X \in L$ to the vector with
$1$ in the $Y$-component if $Y \geq X$ and $0$ otherwise.
The $X$-component of this surjection is
the map $\chi_X: kS \to k$ defined on the elements $y \in S$ by
\begin{gather*}
 \chi_X(y) = \begin{cases}
  1, & \text{ if } \supp(y) \leq X, \\
  0, & \text{ otherwise}.
 \end{cases}
\end{gather*}
The elements
\begin{equation}\label{LRB: idempotents in kL}
E_X = \sum_{Y \geq X} \mu(X,Y) Y 
\end{equation}
in $kL$, one for each $X \in L$, correspond to the standard basis vectors of
$\prod_{X\in L}k$ under the isomorphism $kL \cong \prod_{X \in L} k$.
In the above $\mu$ denotes the M\"obius function of the lattice $L$
\cite[\S3.7]{Stanley1997}. The elements $\{E_X\}_{X \in L}$
form a basis of $kL$ and a \emph{complete system of primitive
orthogonal idempotents} for $kL$ (see the next section for the definition).

\section{Primitive Idempotents of the Semigroup Algebra}
\label{LRB: section: primitive idempotents}

Let $A$ be a $k$-algebra. An element $e \in A$ is \emph{idempotent} if $e^2
= e$.  It is a \emph{primitive idempotent} if $e$ is idempotent and we
cannot write $e = e_1 + e_2$ where $e_1$ and $e_2$ nonzero idempotents in
$A$ with $e_1e_2 = 0 = e_2e_1$.  Equivalently, $e$ is primitive iff $A e$
is an indecomposable $A$-module. A set of elements $\{e_i\}_{i\in I}
\subset A$ is a \emph{complete system of primitive orthogonal idempotents}
for $A$ if $e_i$ is a primitive idempotent for every $i$, if $e_ie_j = 0$
for $i \neq j$ and if $\sum_i e_i = 1$. If $\{e_i\}_{i\in I}$ is a complete
system of primitive orthogonal idempotents for $A$, then $A \cong
\bigoplus_{i\in I} Ae_i$ as left $A$-modules and $A \cong \bigoplus_{i,j\in
I} e_iAe_j$ as $k$-vector spaces.

Let $S$ denote a left regular band with identity. 
For each $X \in L$, fix an $x \in S$ with $\supp(x) =
X$ and define elements in $kS$ recursively by the formula,
\begin{gather}
 \label{LRB: equation: idempotents}
 e_X = x - \sum_{Y > X} x e_Y.
\end{gather}

\begin{Lemma} \label{LRB: idempotent lemma}
 Let $w \in S$ and $X \in L$. If $\supp(w) \not\leq X$, then $w e_X = 0$.
\end{Lemma}
\begin{proof}
We proceed by induction on $X$. This is vacuously true if $X = \hat 1$.
Suppose the result holds for all $Y \in L$ with $Y >
X$. Suppose $w \in S$ and $W = \supp(w) \not\leq X$. Using
the definition of $e_X$ and the identity $wxw = wx$
(LRB2),
\begin{align*}
w e_X & = wx - \sum_{Y > X} wx e_Y = wx - \sum_{Y > X} wx(we_Y).
\end{align*}
By induction, $we_Y = 0$ if $W \not\leq Y$. Therefore, the
summation runs over $Y$ with $W \leq Y$. 
But $Y > X$ and $Y \geq W$ iff $Y \geq W \vee X$, so the summation runs
over $Y$ with $Y \geq W \vee X$. 
\begin{align*}
we_X = wx - \sum_{Y > X} wx(we_Y) = wx - \sum_{Y \geq X \vee W} wx e_Y.
\end{align*}
Now let $z$ be the element of support $X \vee W$ chosen in defining $e_{X
\vee W}$. So $e_{X\vee W} = z - \sum_{Y>X\vee W} ze_Y$. 
Note that $ze_{X\vee W} = e_{X\vee W}$ since $z = z^2$. Therefore, $z =
\sum_{Y\geq X\vee W} ze_Y$. Since $\supp(wx) = W \vee X = \supp(z)$, it
follows from 
Proposition \ref{LRB: associated lattice} (\ref{LRB: xy=x iff XleqY})
that
$wx = wxz$.
Combining the last two statements,
\begin{align*}
we_X = wx - \sum_{Y \geq X \vee W} wx e_Y 
& = wx\left(z - \sum_{Y \geq X \vee W} z e_Y\right) = 0.
\qedhere
\end{align*}
\end{proof}

\begin{Theorem}
 \label{LRB: complete system of primitive orthogonal idempotents}
Let $S$ denote a finite left regular band with identity and $L$ its
support lattice. Let $k$ denote an arbitrary field.
The elements $\{e_X\}_{X \in L}$ form a complete
system of primitive orthogonal idempotents in the semigroup algebra
$kS$.
\end{Theorem}
\begin{proof}
\emph{Complete}. $1$ is the only element of support
$\hat 0$. Hence, $e_{\hat 0} = 1 - \sum_{Y > \hat 0} e_Y$. 
Equivalently, 
$\sum_X e_X = 1$.

\emph{Idempotent.} 
Since $e_Y$ is a linear combination of elements of support at least $Y$, $e_Y
z = e_Y$ for any $z$ with $\supp(z) \leq Y$
(Proposition \ref{LRB: associated lattice} (\ref{LRB: xy=x iff XleqY})).
Using the definition of $e_X$,
the facts $e_X = xe_X$ and $e_Y = e_Yy$, and Lemma \ref{LRB: idempotent lemma},
\begin{align*}
e_X^2 = \left(x - \sum_{Y > X} xe_Y\right)e_X =
xe_X - \sum_{Y > X} xe_Y(ye_X) = xe_X = e_X.
\end{align*}

\emph{Orthogonal.} We show that for every $X \in L$,
$e_Xe_Y = 0$ for $Y \neq X$. If $X = \hat 1$, then $e_X e_Y = e_X x e_Y = 0$
for every $Y \neq X$ by Lemma \ref{LRB: idempotent lemma} since $X = \hat 1$
implies $X \not\leq Y$.  Now suppose the result holds for $Z > X$. That is,
$e_Ze_Y = 0$ for all $Y \neq Z$. If $X \not\leq Y$, then $e_Xe_Y = 0$ by
Lemma \ref{LRB: idempotent lemma}. If $X < Y$, then $e_X e_Y = xe_Y - \sum_{Z > X}
x (e_Ze_Y) = xe_Y - xe_Y^2 = 0.$

\emph{Primitive.} We'll show that $e_X$ lifts $E_X = \sum_{Y\geq
X} \mu(X,Y)Y$ (see Equation (\ref{LRB: idempotents in kL})) for all $X \in L$,
a primitive idempotent in
$kL$. (Then since $e_X$ lifts a primitive idempotent, it is itself a
primitive idempotent.)
If $X = \hat 1$, then $\supp(e_{\hat 1}) = \hat 1 =
E_{\hat 1}$. Suppose the result holds for $Y > X$. Then
$\supp(e_X) = \supp(x - \sum_{Y>X}xe_Y) = X -
\sum_{Y>X}(X\vee E_Y)$. Since $E_Y$ is a linear combination
of elements $Z \geq Y$, it follows that $X \vee E_Y = E_Y$ if $Y>X$.
Therefore, $\supp(e_X) = X - \sum_{Y>X}E_Y$. The M\"obius
inversion formula 
\cite[\S3]{Stanley1997}
applied to $E_X = \sum_{Y\geq X} \mu(X,Y)Y$
gives $X = \sum_{Y \geq X} E_X$. Hence,
$\supp(e_X) = X - \sum_{Y>X}E_Y = E_X$.
%
\end{proof}

\begin{Remark}
 We can replace $x \in S$ in Equation
 (\ref{LRB: equation: idempotents}) with any linear
 combination  $\tilde x = \sum_{\supp(x)=X} \lambda_x x$ 
 of elements of support $X$ whose
 coefficients $\lambda_x$ sum to $1$.
The proofs still hold since the element $\tilde x$ 
is idempotent and satisfies
$\supp(\tilde x) = X$ and $\tilde x y = \tilde x$
for all $y$ with $\supp(y) \leq X$. 
Unless explicitly stated we will use the idempotents constructed
above.
\end{Remark}

\begin{Corollary} \label{LRB: basis of idempotents}
The set $\{ xe_{\supp(x)} \mid x \in S\}$ is a basis of $kS$
of primitive idempotents (not necessarily orthogonal idempotents).
\end{Corollary}
\begin{proof}
Let $y \in S$. Then by Theorem \ref{LRB: complete system of primitive
orthogonal idempotents} and Lemma \ref{LRB: idempotent lemma},
$$
y = y1 = y \sum_Z e_Z = \sum_{Z \geq \supp(y)} ye_Z = \sum_{Z \geq \supp(y)}
(yz)e_Z,
$$
where $z \in S$ was the element used to define $e_Z$.
Since $\supp(yz) = \supp(y) \vee \supp(z) = Z$, every element $y \in S$ 
is a linear
combination of elements of the form $xe_{\supp(x)}$. So the elements
$x e_{\supp(x)}$, one for each $x$ in $S$,
span $kS$. Since the number of these elements is the cardinality of $S$,
which is the dimension of $kS$, the set forms a basis of $kS$. The
elements are idempotent since $(x e_X)^2 = (x e_X)(x e_X) = x e_X^2 = x e_X$
(since $xyx=xy$ for all $x,y\in S$). Since $xe_X$ lifts the
primitive idempotent $E_X = \sum_{Y\geq X} \mu(X,Y)Y \in kL$,
it is also a primitive idempotent (see the end of the proof of
Theorem \ref{LRB: complete system of primitive orthogonal idempotents}).
\end{proof}

\section{Projective Indecomposable Modules of the Semigroup Algebra}
 \label{LRB: section: projective indecomposable modules}

For $X \in L$, let $S_X \subset S$ denote the
set of elements of $S$ of support $X$. For $y \in S$ and $x \in S_X$
define
\begin{gather*}
y\cdot x = \begin{cases}
yx, & \supp(y) \leq \supp(x), \\
0,& \supp(y) \not\leq \supp(x).
\end{cases}
\end{gather*}
Then $\cdot$ defines an action of $kS$ on the $k$-vector space $kS_X$ spanned
by $S_X$. 

\begin{Lemma}
 \label{LRB: basis for indecomposable modules}
Let $X \in L$. Then $\{x e_X \mid \supp(x) = X\}$ is a basis
for $(kS)e_X$.
\end{Lemma}
\begin{proof}
Suppose $\sum_{w \in S} \lambda_w we_X \in kS e_X$. 
If $\supp(w) \not\leq X$, then $we_X = 0$. So suppose $\supp(w) 
\leq X$. Then $\supp(wx) = \supp(w) \vee X = X$.
Therefore,
$$
\sum_{w \in S} \lambda_w we_X = \sum_{w \in S} 
\lambda_w (wx) e_X \in \operatorname{span}_k\{y e_X \mid \supp(y) = X\},
$$
where $x$ is the element chosen in the construction of $e_X$ (recall
that $e_X = xe_X$ since $x^2 = x$).
So the elements span $kS e_X$. These elements are linearly independent
being a subset of a basis of $kS$ (Corollary \ref{LRB: basis of idempotents}).
\end{proof}

\begin{Proposition}
 \label{LRB: projective indecomposable modules}
There is a $kS$-module isomorphism $kS_X \cong kS e_X$ given by
right multiplication by $e_X$.
Therefore, the $kS$-modules $kS_X$ are all the
projective indecomposable $kS$-modules. The radical
of $kS_X$ is $\operatorname{span}_k\{y - y' \mid y,y' \in S_X\}$.
\end{Proposition}
\begin{proof}
Define a map $\phi: kS_X \to kS e_X$ by $w \mapsto
we_X$.  Then $\phi$ is surjective since $\phi(y) = ye_X$
for $y \in S_X$ and since $\{ye_X \mid \supp(y) = X\}$
is basis for 
$kS e_X$ (Lemma \ref{LRB: basis for indecomposable modules}).
Since $\dim kS_X = \#S_X = \dim kS e_X$,
the map $\phi$ is an isomorphism of $k$-vector spaces.

\emph{$\phi$ is a $kS$-module map}. Let $y \in S$
and let $x \in S_X$. If $\supp(y) \leq X$, then $\phi(y \cdot
x) = \phi(yx) = yxe_X = y \phi(x)$. If $\supp(y) \not\leq X$, then $y
\cdot x = 0$. Hence, $\phi(x \cdot y) = 0$. Also, since $\supp(y)
\not\leq X$, it follows 
from Lemma \ref{LRB: idempotent lemma}
that $y e_X = 0$. Therefore, $y\phi(x) = y xe_X = yx(y
e_X) = yx 0 = 0.$ So $\phi(y\cdot x) = y \phi(x)$. Hence $\phi$ is an
isomorphism of $kS$-modules. 

Since all the projective indecomposable
$kS$-modules (upto isomorphism) are of the form $kS e_X$ for a complete
system of primitive orthogonal idempotents $\{e_X\}$, the $kS$-modules $kS_X$
are all the indecomposable projective $kS$-modules.
\end{proof}

\section{The Quiver of the Semigroup Algebra}
\label{LRB: quiver}

Let $A$ be a finite dimensional $k$-algebra whose simple modules are all
$1$-dimensional. The \emph{$\Ext$-quiver} or \emph{quiver} of $A$ is the
directed graph $Q$ with one vertex for each isomorphism class of simple
modules and $\dim_k(\Ext^1_A(M_X,M_Y))$ arrows from $X$ to $Y$, where $M_X$
and $M_Y$ are simple modules of the isomorphism classes corresponding to
the vertices $X$ and $Y$, respectively. The \emph{path algebra} $kQ$ of $Q$
is the $k$-algebra spanned by paths of $Q$ with multiplication induced by
path composition: if two paths in $Q$ compose to form another path, then
that is the product; if the paths do not compose, then the product is 0.
If $Q$ is the quiver of $A$, then there exists a $k$-algebra surjection
from $kQ$ onto $A$. Although the quiver $Q$ is canonical, this surjection
is not.

Let $S$ be a left regular band with identity and let $L$ denote the support
lattice of $S$.
Let $X, Y \in L$ with $Y \leq X$ and fix $y \in S$ with $\supp(y) = Y$. 
Define a relation on the elements of
$S_X$ by $x \smile x'$ if there exists an element $w \in S$ satisfying $y <
w$, $w < yx$ and $w < yx'$. (Equivalently, $yw = w$, $wx = yx$, $wx' =
yx'$ and $\supp(w) < X$.) 
Note that $x \smile x'$ iff $x \smile yx'$. 
Also note that for $X = \hat1$ and $Y = \hat0$, the relation 
becomes
$x \smile x'$ iff there exists $w \neq 1$ such that
$x > w$ and $x' > w$.

The relation $\smile$ is symmetric and reflexive, but not necessarily
transitive. Let $\sim$ denote the transitive closure of $\smile$.
%
Let $a_{XY} = {}^\#(S_X/\sim) - 1$,
the number of equivalence classes of $\sim$ minus one.
If $Y \not\leq X$, define $a_{XY} = 0$. 
In order to avoid confusion, we denote by $a_{XY}^S$ the number $a_{XY}$
computed in $S$.
%
Since $u < v$ implies $yu < yv$ for all $u, v, y \in S$ (follows from
(LRB2)), it follows that the relations $\smile$ and $\sim$ do not depend on
the choice of $y$ with $\supp(y) = Y$.  


\begin{Lemma}
\label{LRB: dimensions of Ext spaces}
Let $S$ be a finite left regular band with identity and $L$ its support 
lattice. Let $M_X$ and $M_Y$
denote the simple modules with irreducible characters $\chi_X$ and
$\chi_Y$, respectively. Then
\begin{align*}
\dim(\Ext^1_A(M_X,M_Y)) = a_{XY}.
\end{align*}
\end{Lemma}
\begin{proof}
The proof is rather lengthy so we postpone it until \S\ref{appendix:
proof of proposition}.
\end{proof}

\begin{Theorem}
Let $S$ be a left regular band with identity and $L$ the support lattice of
$S$. Let $k$ denote a field. The quiver of the semigroup algebra $kS$ has
$L$ as the vertex set and $a_{XY}$ arrows from the vertex $X$ to the vertex $Y$.
\end{Theorem}

\section{An Inductive Construction of the Quiver}
\label{section: an inductive construction of the quiver}

In this section we describe how knowledge about the numbers
$a_{\hat1\hat0}^{S'}$ for certain subsemigroups $S'$ of $S$ determine all
the numbers $a_{XY}^{S}$. This allows for an inductive construction of
the quiver of a left regular band.

Suppose $S$ is a left regular band with identity. Let $X, Y \in L$ with $Y
\leq X$ and let $y \in S$ be an element with $\supp(y) = Y$. Then $yS = \{
yw : w \in S \}$ and $S_{\leq X} = \{ w \in S : \supp(w) \leq X \}$ are 
subsemigroups of $S$.

\begin{Proposition}
Let $S$ be a left regular band with identity, and let $L$ denote the
support lattice of $S$. Suppose $y \in S$ and $X \in L$. The quiver
of the semigroup algebra $k(yS_{\leq X})$ of the left regular band $yS_{\leq
X}$ is the full subquiver of the quiver of the semigroup algebra $kS$ on
the vertices in the interval $[\supp(y), X] \subset L$.
\end{Proposition}

The Proposition follows from the following Lemma that shows the 
number of arrows from $X$ to $Y$ in the quiver of $kS$ is the number
of arrows from $\hat 1$ to $\hat 0$ in the quiver of $k(yS_{\leq X})$,
where $y \in S$ is any element of support $Y$.
Recall that
$a_{\hat1\hat0}^{yS_{\leq X}}$ denotes the number $a_{\hat1\hat0}$ computed in
the left regular band $yS_{\leq X}$.
\begin{Lemma}
Let $S$ be a left regular band with identity. Then $a^S_{XY} =
a_{\hat1\hat0}^{yS_{\leq X}}$. 
That is, the number $a_{XY}$ computed in $S$ is the number $a_{\hat1\hat0}$ 
computed in $yS_{\leq X}$.
\end{Lemma}
\begin{proof}
If $\supp(y) \not\leq X$, then $yS_{\leq X}$ is empty. So $a_{X,Y}^S = 0 =
a_{\hat1\hat0}^{yS_{\leq X}}$. So suppose $\supp(y) \leq X$.

Since $x \sim x'$ iff $x \sim yx'$ for any elements $x, x'$ of support $X$,
every equivalence class of $\sim$ (on $S_X$) contains an element of $yS_X$.
Therefore, $a_{XY} + 1$ is the number of equivalence classes of $\sim$
restricted to $yS_X$.

Since $yS_{\leq X}$ is a subsemigroup of $S$, the support lattice of
$yS_{\leq X}$ is the image of $yS_{\leq X}$ in $L$. Therefore,
the support lattice of $yS_{\leq X}$ is the interval $[Y,X]$ in $L$.
Since the top and bottom elements of $[Y,X]$ are $X$ and $Y$
respectively, the number $a_{\hat1\hat0}^{yS_{\leq X}} + 1$ is the number
of equivalence classes of $\sim$ restricted to $yS_X$.
\end{proof}

Therefore, if the numbers $a_{\hat1\hat0}^{yS_{\leq X}}$ are known for all the
subsemigroups of $S$ of the form $yS_{\leq X}$, then the quiver of
$kS$ is known. We illustrate this technique with two examples in
the next two sections.

\section{Example: The Free Left Regular Band}
\label{section: example: the free left regular band}

Let $S = F(A)$ denote the free left regular band on a finite set $A$ (defined
in Example \ref{LRB: free LRB}). Recall that the support lattice $L$ of $S$
is the set of subsets of $A$. 

Let $y \in S$ and let $Y \subset A$ denote the set of elements occuring in
the sequence $y$. Then $yS$ is the set of all sequences of elements of
$A$ (without repetition) that begin with the sequence $y$.  Therefore, $yS$
is isomorphic to the free left regular band on $A \backslash Y$.  If $X
\subset A$ (so $X \in L$), then $S_{\leq X}$ is the set of all sequences
containing only elements from $X$ (without repetition). Therefore, $S_{\leq
X}$ is also a free left regular band.  It follows that $yS_{\leq X}$ is a
free left regular band for any $y \in S$ and $X \subset A$. Therefore, the
quiver of $S$ is determined once the numbers $a_{\hat0\hat1} =
a_{A\emptyset}$ are known for any free left regular band. 

If two sequences $x, y \in S$ begin with the same element $a \in A$, then
$ax = x$ and $ay = y$. Therefore, $x \sim y$.  Conversely, if $x \smile y$,
then there is a nonempty sequence $w$ such that $wx = x$ and $wy = y$. Then
$x$ and $y$ both begin with the first element of $w$. Therefore, $x \sim y$
iff $x$ and $y$ are sequences begining with the same element.  So the
equivalence classes of $\sim$ are determined by the first elements of the
sequences in $S$. Hence, $a_{\hat1\hat0} = {}^\#(A) - 1$. This argument
applies to any free left regular band with identity, so $a_{XY} = {}^\#(X
\backslash Y) - 1$ since $yS_{\leq X}$ is isomorphic to the free left
regular band on the elements $X \backslash Y$.

\begin{Theorem}\emph{(K. S. Brown, private communication.)}
Let $S=F(A)$ be the free left regular band on a finite set $A$ and let $k$
denote a field. Then the quiver of the semigroup algebra $kS$ has one
vertex $X$ for each subset $X$ of $A$ and ${}^\#(X \backslash Y) - 1$
arrows from $X$ to $Y$ if $Y \subset X$ (and no other arrows or vertices).
\end{Theorem}

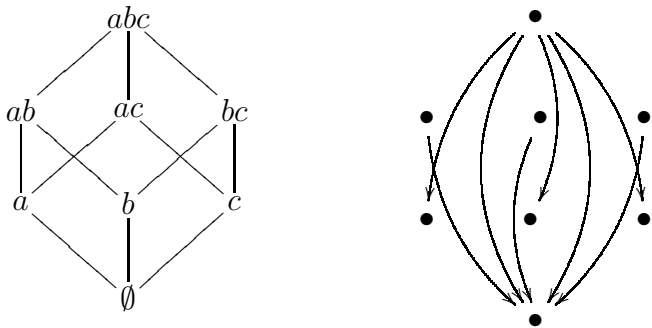
\begin{figure}
\begin{gather*}
\xymatrix@M=1.5pt{
  & abc \ar@{-}[dl] \ar@{-}[rd] \ar@{-}[d] \\
ab \ar@{-}[d]\ar@{-}[dr] & ac \ar@{-}[dl]\ar@{-}[dr] & bc \ar@{-}[d]\ar@{-}[dl] \\
a \ar@{-}[dr]  & b \ar@{-}[d]  & c \ar@{-}[dl] \\
 & \emptyset
}
\hspace{5em} 
\xymatrix{
  & \bullet \ar@/^4ex/[ddd] \ar@/_4ex/[ddd] \ar@/_2ex/[ddl] \ar@/^2ex/[rdd]
\ar@<-0.4ex>@/^2ex/[dd] \\
\bullet \ar@/_2ex/[ddr] & \ \bullet \ar@<0.4ex>@/_2ex/[dd] & \bullet \ar@/^2ex/[ddl] \\
\bullet & \bullet \ & \bullet \\
 & \bullet
}
\end{gather*}
\caption[The quiver of the free left regular band on three generators]
{
The support lattice and the quiver of the semigroup algebra of the free left regular band on three generators. See Figure \ref{figure: free lrb}.}
\end{figure}

\section{Example: The Face Semigroup of a Hyperplane Arrangement}
\label{section: example: the face semigroup of a hyperplane arrangement}

Let $\FF$ be the face semigroup of a central hyperplane arrangement $\AA$
and let $\LL$ be the intersection lattice of $\AA$
(see Example \ref{LRB: hyperplane arrangements}.)
Let $X, Y \in \LL$ and $y$ a face of support $Y$. Then the subsemigroup
$y\FF_{\leq X}$ is the semigroup of faces of a hyperplane arrangement with
intersection lattice $[Y,X] \subset \LL$. (Specifically, this hyperplane
arrangement is given by $\{X \cap H : H \in \AA, Y \subset H, X \not\subset
H\}$.) Therefore, we know all the numbers $a_{XY}$ for $\FF$ if we know the
number $a_{\hat1\hat0}$ for the face semigroup of an arbitrary arrangement.

If $\LL$ contains only one element, then $\hat0 = \hat1$ and
$a_{\hat1\hat0} = 0$. Suppose that $\LL$ contains at least two elements.
It is well-known that for any two distinct chambers $c$ and $d$, there
exists a sequence of chambers $c_0 = c, c_1, \ldots, c_i = d$ such that
$c_{j-1}$ and $c_j$ share a common codimension one face $w_j$ for each $1
\leq j \leq i$ \cite[\S I.4E Proposition 3]{Brown1989}. 
Therefore, $c_{j-1} \smile c_j$ unless $w_j$ is of support
$\hat0$, in which case $\LL$ has two elements. 
Equivalently, $c \sim d$ iff the arrangement is of rank
greater than 2. So if $\LL$ has exactly two elements, then
$a_{\hat1\hat0} = 1$ and if $\LL$ has more than two elements then
$a_{\hat1\hat0} = 1$.

\begin{Theorem}[{\cite[Corollary 8.4]{Saliola2006:FaceSemigroupAlgebra-arxiv}}]
The quiver $Q$ of the semigroup algebra $k\FF$ coincides with the Hasse diagram
of $\LL$. That is, there is exactly one arrow $X \to Y$ iff $Y \lessdot X$.
\end{Theorem}

In \cite{Saliola2006:FaceSemigroupAlgebra-arxiv} the relations of the
quiver are also determined. Let $I$ be the ideal generated by the
following elements, one for each interval $[Z,X]$ of length two in $\LL$,
\begin{gather*}
 \sum_{Y: Z \lessdot Y \lessdot X} X \to Y \to Z. 
\end{gather*}
Then $k\FF \cong kQ/I$ as $k$-algebras, where $kQ$ is the path algebra of
$Q$.

\section{Idempotents in the subalgebras $k(yS)$ and $kS_{\geq X}$}
\label{section: idempotents in the subalgebras}

This section describes the subalgebras of $kS$
generated by the subsemigroups $yS$ and $S_{\leq Y}$ of $S$.

Let $S$ be a left regular band. Recall that for $y \in S$, the set $yS = \{
yw : w \in S \} = \{w \in S : w > y \}$ is a subsemigroup of $S$ (and hence
a left regular band).  Note that if $\supp(y') = \supp(y)$ then the left
regular bands $yS$ and $y'S$ are isomorphic with isomorphism given by
multiplication by $y$ (the inverse is multiplication by $y'$).  Since $yS$
is a subsemigroup of $S$, the support lattice of $yS$ is the image of $yS$
in $L$ by Propositon \ref{LRB: associated lattice}, which is the interval
$[Y,\hat 1]$. 

\begin{Proposition}
Let $S$ be a left regular band, let $y \in S$ and let $Y = \supp(y)$.
There exists a complete system of primitive orthogonal idempotents $\{e_X
: X \in L\}$ in $kS$ such that $\{e_X : X \geq Y\}$ is a complete system of
primitive orthogonal idempotents in the semigroup algebra $k(yS)$.
Moreover, $k(yS) = (\sum_{X \geq Y} e_X) kS$.
\end{Proposition}
\begin{proof}
For each $X \in L$, fix $x \in S$ with $\supp(x) = X$.  If $X \geq Y$, then
replace $x$ with $yx$. Note that $\supp(yx) = \supp(x)$ since $X \geq Y$.
Therefore, $x > y$ if $X \geq Y$. The formula $e_X = x - \sum_{W > X}
xe_W$ for $X \in L$ defines a complete system of primitive orthogonal
idempotents for $kS$ (Theorem \ref{LRB: complete system of
primitive orthogonal idempotents}).
And since the support lattice of $yS$ is $[Y,\hat1]
\subset L$, the elements
$e_X = x - \sum_{W > X} xe_W$ for $X \geq Y$ define a complete system of
primitive orthogonal idempotents in $k(yS)$.
Since $y$ is the identity of $yS$, we have
$y = \sum_{X \geq Y} e_X$. Therefore,
$k(yS) = y(kS) = (\sum_{X \geq Y} e_X) kS$.
\end{proof}

If $Y \in L$, then $S_{\leq Y} = \{ w \in S : \supp(w) \leq Y \}$ is a
subsemigroup of $S$. The support lattice of $S_{\leq Y}$ is
the interval $[\hat 0, Y]$ of $L$. Let $\operatorname{proj}_{kS_{\leq
Y}}: kS \to kS_{\leq X}$ denote the projection onto the subspace
$kS_{\leq X}$ of $kS$.

\begin{Proposition}
Let $S$ be a left regular band and $Y \in L$. 
Let $\{e_X : X \in L\}$ denote a complete system of primitive orthogonal
idempotents of $kS$. Then $\{\operatorname{proj}_{kS_{\leq Y}}(e_X) : X \leq
Y\}$ is a complete system of primitive orthogonal idempotents of $kS_{\leq
Y}$. Moreover, the semigroup algebra
$k(S_{\leq Y})$ is isomorphic to $kS(\sum_{X \leq Y} e_X)$.
\end{Proposition}
\begin{proof}
The map $\operatorname{proj}_{kS_{\leq Y}}$ is an algebra morphism $kS \to
kS_{\leq Y}$. This follows from the fact that $\supp(wx) = \supp(w) \vee
\supp(x)$ for any $x, w \in S$.  So if $X \leq Y$, then
$\operatorname{proj}_{kS_{\leq Y}}(e_X) = x - \sum_{W > X} x
\operatorname{proj}_{kS_{\leq Y}}(e_W)$ since $e_X = x - \sum_{W > X}
xe_W$.  Therefore, the elements $\operatorname{proj}_{kS_{\leq Y}}(e_X)$
for $X \leq Y$ form a complete system of primitive orthogonal idempotents
for the semigroup algebra of the left regular band $S_{\leq Y}$ (Theorem
\ref{LRB: complete system of primitive orthogonal idempotents}).
Since $\operatorname{proj}_{kS_{\leq Y}}$ is an algebra morphism, it
restricts to a surjective morphism of algebras
$\operatorname{proj}_{kS_{\leq Y}}: kS(\sum_{X \leq Y} e_X) \to k(S_{\leq
Y})$. Since $kS_X \cong (kS)e_X$ for all $X \in L$ as $kS$-modules
(Proposition \ref{LRB: projective indecomposable modules}),
$\dim(kS_{\leq Y}) = \dim(\sum_{X \leq Y} (kS)e_X)$. So
$\operatorname{proj}_{kS_{\leq Y}}$ is an isomorphism. Its inverse
is right multiplication by $\sum_{X \leq Y} e_X$.
\end{proof}

\section{Cartan Invariants of the Semigroup Algebra}
\label{LRB: section: Cartan invariants}

The \emph{Cartan invariants} of a finite dimensional
$k$-algebra $A$ are the numbers $\dim_k(\Hom_{A}(A e_X, A e_Y))$, where
$\{e_X\}_{X \in I}$ is a complete system of primitive orthogonal
idempotents for $A$. They are independent of the choice of $\{e_X\}_{X \in
I}$.

Let $S$ be a left regular band with identity and let $L$ denote the support
lattice of $S$. For $X,Y \in L$, define numbers $m(Y,X)$ 
follows. If $Y \not\leq X$, then $m(Y,X) = 0$. 
If $Y \leq X$, then define $m(Y,X)$ by the formulas
\begin{equation}
\label{Equation: definition of m(Y,X)}
\sum_{W \leq Y \leq X} m(Y,X) = {}^\#(wS_X),
\end{equation}
one for each $W \in L$, where $w$ is an element of support $W$. (Recall that
the number ${}^\#(wS_X)$ does not depend on the choice of $w$ with
$\supp(w) = W$.) Equivalently, 
$$m(Y,X) = \sum_{Y \leq W \leq X} \mu(Y,W) \ {}^\#(wS_X),$$
where $\mu$ is the M\"obius function of $L$ 
\cite[\S3.7]{Stanley1997}.

\begin{Proposition}
 \label{LRB: cartan invariants}
Let $S$ be a left regular band with identity. Let $\{e_X\}_{X \in L}$
denote a complete system of primitive orthogonal idempotents for $kS$. 
Then for any $X, Y$,
 $$
 \dim(e_Y kS e_X) = \dim \Hom_{kS}(kS e_Y, kS e_X) = m(Y,X).
 $$
Therefore, the numbers $m(Y,X)$ are the Cartan invariants of $kS$.
\end{Proposition}
\begin{proof}
The first equality follows from the identity $\Hom_A(Ae, Af) \cong eAf$ for
idempotents $e,f$ of a $k$-algebra $A$. 
If $Y \not\leq X$, then it follows 
from (LRB2) and Lemma \ref{LRB: idempotent lemma}
that $e_Y kS e_X = 0$. Suppose that $Y \leq X$.
From the previous section,
$k(yS) = \sum_{W \geq Y} e_W kS$ for some complete system of primitive
orthogonal idempotents. Combined with the isomorphism $kS_X \cong kSe_X$
we get 
$k(yS_X) \cong \bigoplus_{Y \leq W \leq X} e_W kS e_X$.
Therefore,
$$
\sum_{Y \leq W \leq X} m(W,X)
= \dim(k(yS_X)) = \sum_{Y \leq W \leq X} \dim(e_W kS e_X).
$$
The result now follows by induction. If $X = Y$,
then $\dim e_X kS e_X = m(X,X)$. Suppose the result
holds for all $W$ with $Y < W \leq X$. Then
\begin{align*}
 \dim e_Y kS e_X
 &= \sum_{Y \leq W \leq X} m(W,X) - \sum_{Y < W \leq X} \dim e_W kS e_X \\
 &= \sum_{Y \leq W \leq X} m(W,X) - \sum_{Y < W \leq X} m(W,X) \\
 &= m(Y,X). 
\qedhere
\end{align*}
\end{proof}

\section{Example: The Face Semigroup of a Hyperplane Arrangement}
\label{section: the cartan invariants for hyperplane arrangements}
Let $\FF$ denote the semigroup of faces of a hyperplane arrangement $\AA$.
Then $ {}^\#(w\FF_X)$ is the number of faces of support $X$ containing $w$
as a face. Zaslavsky's Theorem \cite{Zaslavsky1975} gives that this is
$\sum_{W \leq Y \leq X} |\mu(Y,X)|$, where $\mu$ is the M\"obius function
of the intersection lattice of $\AA$. Comparing this with Equation
(\ref{Equation: definition of m(Y,X)}) we conclude that the Cartan
invariants of $k\FF$ are $m(Y,X) = |\mu(Y,X)|$. These were also computed in
\cite[Proposition 6.4]{Saliola2006:FaceSemigroupAlgebra-arxiv}.

\section{Example: The Free Left Regular Band}
\label{section: the cartan invariants for a free left regular band}
Let $S$ be a free left regular band on a finite set $A$. The support
lattice of $S$ is the lattice of subsets of $A$. Therefore, $\mu(Y,W) =
(-1)^{{}^\#(W\backslash Y)}$ \cite[Example 3.8.3]{Stanley1997} for
any $Y, W \in L$. And ${}^\#(wS_X) = {}^\#(X\backslash W)!$ since
the number of elements of maximal support in the free left regular band
on $A$ is precisely ${}^\#A!$. 
If $n = {}^\#X$ and $j = {}^\#Y$, and $Y \subset X$, then 
\begin{align*}
m(Y,X)
& = \sum_{Y \leq W \leq X} \mu(Y, W) \ {}^\#(wS_X) \\
& = \sum_{Y \leq W \leq X} (-1)^{{}^\#W - j} \ (n - {}^\#W)!\\
& = \sum_{i = j}^{n} \sum_{Y \subset W \subset X \atop {}^\#W = i}
(-1)^{i - j} (n - i)! \\
& = \sum_{i = j}^{n} (-1)^{i-j} (n - i)! \binom{n-j}{i - j} \\
& = (n-j)! \sum_{i = j}^{n} \frac{(-1)^{i-j}}{(i-j)!} \\
& = (n-j)! \sum_{i = 0}^{n - j} \frac{(-1)^{i}}{i!}.
\end{align*}
Therefore, the number $m(Y,X)$ depends only on the cardinality of $X
\backslash Y$ and we denote it by $m_i$ where $i = {}^\#(X \backslash Y)$.

We will now prove that these numbers count paths in the quiver of $kS$.
For a set $A$ of cardinality $n$, let $Q_n$ be the directed graph with one
vertex for each subset of $A$ and ${}^\#(X\backslash Y)-1$ arrows from $X$
to $Y$ if $Y \subset X$. Let $p_n$ denote the number of paths in $Q_n$
beginning at $A$ and ending at $\emptyset$. Note that if $Y \subset X
\subset A$, then the number of paths beginning at $X$ and ending at $Y$ in
$Q_n$ is $p_m$ where $m = {}^\#(X\backslash Y)$.

For each $0 \leq i \leq n-1$ there are $n - i - 1$ arrows from $A$ to sets
of size $i$, and there are $\binom{n}{i}$ such sets, so 
$p_n = \sum_{0 \leq i \leq n-1} \binom{n}{i}(n-i-1) p_i$ 
for $n \geq 1$. Equivalently, 
\begin{align*}
\sum_{0 \leq i \leq n} \binom{n}{i} p_i = 
\sum_{0 \leq i \leq n-1} \binom{n}{i}(n-i) p_i.
\end{align*}
If the $m_i$ satisfy the above recurrence, then $m_i = p_i$
for all $i$ since $m_0 = 1 = p_0$. Well,
\begin{align*}
& \sum_{0 \leq i \leq n-1} \binom{n}{i}(n-i) m_i \\
&= \sum_{0 \leq i \leq n-1} \frac{n!}{(n-i-1)!\ i!} m_i \\
&= \sum_{0 \leq i \leq n-1} \frac{n!}{(n-i-1)!\ i!} 
  \left(i! \sum_{0 \leq j \leq i}\frac{(-1)^j}{j!} \right) \\
&= \sum_{0 \leq i \leq n-1} \frac{n!}{(n-i-1)!} 
  \left(\sum_{0 \leq j \leq i}\frac{(-1)^j}{j!} \right) \\
&= \sum_{1 \leq k \leq n} \frac{n!}{(n-k)!} 
  \left(\sum_{0 \leq j \leq k-1}\frac{(-1)^j}{j!} \right) \\
&= \sum_{1 \leq k \leq n} \frac{n!}{(n-k)!} 
  \left(\sum_{0 \leq j \leq k}\frac{(-1)^j}{j!} - \frac{(-1)^k}{k!}\right) \\
&= \sum_{1 \leq k \leq n} \frac{n!}{(n-k)!} 
  \left(\sum_{0 \leq j \leq k}\frac{(-1)^j}{j!}\right)
  - \sum_{1 \leq k \leq n} \frac{n!}{(n-k)!} 
  \left(\frac{(-1)^k}{k!}\right) \\
&= \sum_{1 \leq k \leq n} \binom{n}{k} 
  \left(k! \sum_{0 \leq j \leq k}\frac{(-1)^j}{j!}\right)
  - \sum_{1 \leq k \leq n} \binom{n}{k} 
  (-1)^k \\
&= \sum_{1 \leq k \leq n} \binom{n}{k} m_k + 1 \\
&= \sum_{1 \leq k \leq n} \binom{n}{k} m_k + \binom{n}{0} m_0 \\
&= \sum_{0 \leq k \leq n} \binom{n}{k} m_k.
\end{align*}

\begin{Theorem}\emph{(K. S. Brown, private communication.)}
Let $S = F(A)$ be the free left regular band on a finite set $A$.
Then $kS \cong kQ$, where $kQ$ is the path algebra of the quiver
$Q$ of $kS$.
\end{Theorem}
\begin{proof}
Since $Q$ is the quiver of $kS$, there is an algebra surjection $kQ \to
kS$, where $kQ$ is the path algebra of $Q$.  The canonical basis for $kQ$
is the set of paths in $Q$, so using the fact that $m(Y,X) = \dim(e_Y kS
e_X)$ counts the number of paths in $Q$ from $X$ to $Y$ (see the preceeding
two paragraphs), we have $\dim(kQ)
= \sum_{Y, X} m(Y,X) = \sum_{Y,X} \dim( e_Y kS e_X) = \dim(kS)$.
\end{proof}

\section{Future Directions}
\label{section: future directions}

We conclude this paper by providing a few problems for future exploration.

Although this paper successfully determines the quiver of the semigroup
algebra of a left regular band, it says nothing about the quiver relations.
\emph{Describe the quiver relations of the semigroup algebra of a left
regular band with identity.}

The face semigroup algebra of a hyperplane arrangement is a Koszul algebra
\cite[Proposition 9.4]{Saliola2006:FaceSemigroupAlgebra-arxiv} and its
Koszul dual is the incidence algebra of the opposite lattice of the support
lattice of the semigroup. Since this algebra is the semigroup algebra of a
left regular band, it is natural to ask this question for all left regular
bands. \emph{Determine which class of left regular bands give Koszul
semigroup algebras and identify their Koszul duals.} One source of examples
of left regular bands giving Koszul algebras comes from \emph{interval
greedoids} \cite{BjornerZiegler1992:Greedoids}. This will be explored in an
upcoming paper.

Another nice property of the face semigroup algebra of a hyperplane
arrangement is that the quiver of the semigroup algebra coincides with the
support lattice of the semigroup. In fact, the support lattice completely
determines the semigroup algebra. \emph{Determine the left regular bands
$S$ for which the quiver of $kS$ coincides with the support lattice of
$L$.} (From our description of the quiver of $kS$, we have a description of
these left regular bands in terms of the equivalence classes of $\sim$.)
\emph{Determine those $S$ for which the support lattice $L$ completely
determines $kS$.}

A \emph{band} is a semigroup $B$ satisfying $b^2 = b$ for all $b \in B$.
Since left regular bands are bands it is natural to try to generalize these
results to arbitrary bands.  \emph{Describe the quiver of the semigroup
algebra $kB$ of a band $B$ with identity. Construct a complete system of
primitive orthogonal idempotents for $kB$. Determine the bands $B$ for
which $kB$ is a Koszul algebra.}

\section{Appendix: Proof of Lemma \ref{LRB: dimensions of Ext
spaces}}
\label{appendix: proof of proposition}

\begin{LemmaAppendix1}
Let $S$ be a finite left regular band with identity and $L$ its support 
lattice. Let $M_X$ and $M_Y$
denote the simple modules with irreducible characters $\chi_X$ and
$\chi_Y$, respectively. Then
\begin{gather*}
\dim(\Ext^1_A(M_X,M_Y)) = a_{XY}.
\end{gather*}
\end{LemmaAppendix1}
\begin{proof}
As a vector space $M_X = k$ and the action of $kS$ on $M_X$ is
given by $\chi_X$: if $y \in S$ and $\lambda \in k$, then $y \cdot \lambda
= \chi_X(y) \lambda$.

Since the following is a short exact sequence of $kS$-modules
with $kS_X$ projective,
\begin{align*}
0 \longrightarrow 
\ker\left(\chi_X|_{kS}\right) \longrightarrow
kS_X \tot{\chi_X}{\longrightarrow} 
M_X \longrightarrow 0
\end{align*}
Proposition $7.2$ in Chapter V of 
\cite{CartanEilenberg}
gives the exact sequence
\begin{align*}
\Hom_{kS}(kS_X, M_Y) \longrightarrow
\Hom_{kS}&(\ker \left(\chi_X|_{kS}\right), M_Y) \\
& \hspace{3em} \longrightarrow \Ext^1_{kS}(M_X, M_Y) \longrightarrow
0.
\end{align*}
Let $K$ denote the kernel of
$\chi_X|_{kS_X}$. Then $K$ is spanned by the differences of elements of
support $X$. If $f \in \Hom_{kS}(K,\ M_Y)$ and $x, x'$ are elements of
support $X$, then $f(x - x') = 1 f(x - x') = \chi_Y(y) f(x - x') = y \cdot
f(x - x') = f(y \cdot (x - x'))$, for any element $y$ of support $Y$. So if
$Y \not\leq X$ or if $Y = X$, then $f = 0$. Therefore,
$\Hom_{kS}(K,M_Y) = 0$ if $Y \not< X$. It follows that
\begin{align*}
\Ext^1_{kS}(M_X,M_Y) = 0 = a_{XY} \text{ for } Y \not< X.
\end{align*}

Suppose $Y < X$. If $f \in \Hom_{kS}(kS_X,M_Y)$, then 
for all $x \in S_X$, $f(x) = f(x^2)
= f(x \cdot x) = x \cdot f(x) = \chi_Y(x) f(x) = 0 f(x) = 0$ for all $x \in
S$ with $\supp(x) = X$. Therefore, $\Hom_{kS}(kS_X,M_Y)
= 0$. Hence,
\begin{align*}
\Ext^1_{kS}(M_X,M_Y) \cong \Hom_{kS}(K, M_Y) \text{ for } Y < X.
\end{align*}
Suppose $x \smile x'$. Then there exists a $w \in S$ with $y < w$,
$\supp(w) < X$, $wx = yx$ and $wx = yx'$.  Then $x - x' \in K$, and for any
$f \in \Hom_{kS}(K,\ M_Y)$ we have $f(x - x') = \chi_Y(y) f(x - x') = f(yx
- yx') = f(wx - wx') = f(w \cdot (x - x')) = w \cdot f(x - x') = \chi_Y(w)
f(x - x') = 0 f(x - x') = 0$. Therefore, $f(x - x') = 0$ if $x \smile x'$.
If $x \sim x'$, then there exist $x_0 = x, x_1, \ldots, x_i = x'$ such that
$x_{j-1} \smile x_j$ for $1 \leq j \leq i$, and $f(x - x') = f(x_0 - x_1) +
f(x_1 - x_2) + \cdots + f(x_{i-1} + x_i) = 0$. Therefore, $f(x - x') = 0$
if $x \sim x'$. So $f$ can only be nonzero on differences of elements 
in different equivalence classes of $\sim$. 
Moreover, the equivalence classes determine $f$:
if $u \sim x$ and $u' \sim x'$, then 
$f(u - u') = f(u - x) + f(x - x') + f(x' - u') = f(x - x')$.
Therefore, 
\begin{align*}
\dim(\Ext^1_{kS}(M_X,M_Y)) = \dim(\Hom_{kS}(K,M_Y)) \leq a_{XY}.
\end{align*}

Fix $y$ with $\supp(y) = Y$ and let $x, x' \in S_X$
with $x \not\sim x'$. 
Since $\{u - x : u \neq x, \supp(u) = X \}$ is
a basis for $K$, we get a well-defined linear function $f: K \to k$ by
defining
\begin{align*}
f(u - x) =
\begin{cases}
1, & \text{if } u \sim x', \\
0, & \text{otherwise}.
\end{cases}
\end{align*}
We now show that
$f: K \to M_Y$ is a $kS$-module map. That is,
$f(w \cdot (u - x)) = \chi_Y(w) \cdot f(u - x)$ for all
$w \in S$ and for all $u \in S_X$.

Suppose $\supp(w) \not\leq Y$. Then $w \cdot f(u - x) = 0$ since $w$ acts
trivially on $M_Y$. 
If $\supp(w) \not< X$, then $w$ acts trivially on $K$ and so $w \cdot f(u -
x) = 0 = f(w \cdot (u - x))$. 
So suppose $\supp(w) < X$. Then $f(w \cdot (u - x)) = f(wu - wx) = f(wu -
x) - f(wx - x)$. Since $v \sim x'$ iff $yv \sim x'$ for any $v \in S_X$, it
follows that $f(wu - x) = f(ywu - x)$ and $f(wx - x) = f(ywx - x)$. 
If $\supp(yw) = X$, then $ywu = yw = ywx$ (LRB2), so $f(w \cdot (u - x)) = 0$.
If $\supp(yw) < X$, then we have an element $v = yw$ satisfying $v  >
y$, $\supp(v) < X$, $v(wu) = y(wu)$ and $v(wx) = y(wx)$. That is, $wu \sim
wx$ and it follows that $f(wu - x) = f(wx - x)$. So $f(w \cdot (u - x)) =
0$.

Suppose $\supp(w) \leq Y$. Then $w$ acts as the identity on $M_Y$. Hence,
$w \cdot f(u - x) = f(u - x)$. Since $\supp(w) \leq Y$ and $Y \leq X$, we
have that $\supp(w) \leq X$. Therefore, $f(w \cdot (u - x)) = f(wu - wx) =
f(wu - x) - f(wx - x)$. Since $v \sim x'$ iff $yv \sim x'$, we have $f(wu -
x) = f(y(wu) - x) = f(yu - x) = f(u - x)$ since $\supp(w) \leq Y$.
Similarly, $f(wx - x) = f(x - x) = 0$. Therefore,
$f(w \cdot (u - x)) = f(u - x)$.

This establishes that $f: K \to M_Y$ is a $kS$-module map. And since $f$
is nonzero only on differences of the form $u - u'$ with $u \sim x$ and $u'
\sim x'$, there are exactly $a_{XY}$ such $kS$-module maps. These maps are
linearly independent, therefore
\begin{align*}
\dim(\Ext^1_{kS}(M_X,M_Y)) = \dim(\Hom_{kS}(K, M_Y)) \geq a_{XY}. \quad \qedhere
\end{align*}
\end{proof}

\bibliographystyle{hapalike}
\bibliography{references} 

\end{document}